\documentclass[a4paper]{article}

\usepackage[english]{babel}
\usepackage[utf8x]{inputenc}
\usepackage{amsmath}
\usepackage{amsthm}
\usepackage{graphicx}
\usepackage[colorinlistoftodos]{todonotes}

\title{On Elementary Methods To Evaluate Values of The Riemann Zeta Function and Another Closely Related Infinite Series At Natural Numbers}
\author{Dhrushil Badani \\badanidhrushil@gmail.com}

\begin{document}
\maketitle

\begin{abstract}
In this paper, an elementary method to find the values of the Riemann Zeta function at even natural numbers, and to find values of a closely related series at odd natural numbers is presented. Another method, specifically for the evaluation of $\zeta(2),$ is also presented.
\end{abstract}

\section{Introduction}
\begin{quotation}
\emph{''The zeta function  is probably the most challenging and mysterious object of modern mathematics, in spite of its utter simplicity.''}
\flushright
\emph{M.C Gutzwiller.}
\end{quotation}

The Riemann Zeta function is one of the most important functions in Mathematics – it is deeply related to the prime number theorem and has wide-ranging applications in physics, probability theory, analytic number theory and other fields of mathematics. The Riemann Hypothesis, one of the Millennium Prize problems, is closely related to the Riemann Zeta function.

In this paper,the only considered cases are where the argument is purely real. Firstly, we suggest an alternate method to evaluate $\zeta(2)$ i.e. $\frac{\pi^2}{6}$ by setting up a definite integral and evaluating it using methods of calculus. Secondly, we suggest an elementary method by means of which we generalize the procedure to find the sum(s) of the following infinite series:

When $s$ is an even natural number, $$\zeta(s)= \sum_{n=1}^{\infty} {\frac{1}{n^s}}\ $$

and,

When $s$ is an odd natural number, let $$\Phi(s)= \sum_{n=1}^{\infty} {\frac{{(-1)}^n}{{(2n+1)}^s}}\ $$
\pagebreak

\section{An Alternate Method to Calculate $\zeta(2)$}
In this section, we use the tools of integral calculus to set up a definite integral for $\zeta(2)$, and evaluate it.
$$ \zeta(2) = \sum_{n=1}^{\infty} {\frac{1}{n^2}}\ $$

Introduce a parameter $x$ and consider the following function $f(x)$:
$$f(x)=\sum_{n=1}^{\infty} {\frac{x^n}{n^2}}\ $$

Differentiate both sides of this equation with respect to $x$:
$$f'(x)=\sum_{n=1}^{\infty} {\frac{nx^{n-1}}{n^2}}\ $$
$$f'(x)=\frac{1}{x} \sum_{n=1}^{\infty} {\frac{x^n}{n}}\ $$

Also, consider the Taylor Series expansion of $\ln(1+x)$ :
$$\ln(1+x)=-\sum_{n=1}^{\infty} {\frac{{(-1)}^n {x^n} }{n}}\ $$

If we replace $x$ with $-x$ in this expansion, we have,
\begin{equation}
ln(1-x)=-\sum_{n=1}^{\infty} {\frac{x^n}{n}}\
\end{equation}

Thus,
$$ f'(x)=-\frac{\ln(1-x)}{x} $$
$$\frac{ \,\mathrm{d} (f(x))}{ \,\mathrm{d} x}=-\frac{\ln(1-x)}{x}$$
$$\,\mathrm{d}(f(x))=-\frac{\ln(1-x)}{x} \,\mathrm{d}x $$

Integrate both sides of this equation with respect to $x$ from $x=0$ to $x=1$:
$$f(1)-f(0)=-\int_0^1 \frac{\ln(1-x)}{x} \,\mathrm{d}x $$

Now introduce a parameter '$t$' and set
$$I(\alpha)=\int_0^1 \frac{ln{(1-\alpha x+x^2)}}{x} \,\mathrm{d}x$$

Thus,
$$I(2)=\int_0^1 \frac{ln{(1-2x+x^2)}}{x} \,\mathrm{d}x$$\
$$I(2)=\int_0^1 \frac{ln{{(1-x)}^2}}{x} \,\mathrm{d}x$$\
\begin{equation}
I(2)= 2 \int_0^1 \frac{\ln(1-x)}{x} \,\mathrm{d}x
\end{equation}

And,
$$I(-2)= \int_0^1 \frac{ln{(1+2x+x^2)}}{x} \,\mathrm{d}x$$\
$$I(-2)=\int_0^1 \frac{ln{{(1+x)}^2}}{x} \,\mathrm{d}x$$\
\begin{equation}
I(-2)= 2 \int_0^1 \frac{\ln(1+x)}{x} \,\mathrm{d}x 
\end{equation}

Now, set
$$J=\int_0^1 \frac{\ln(1-x)}{x} \,\mathrm{d}x $$

Thus,
\begin{equation}
I(2)=2J
\end{equation}

Let $x=t^2$
Thus,
$$J=\int_0^1 \frac{\ln(1-t^2)}{t^2} \,\mathrm{d}{(t^2)}$$
$$J=2\int_0^1 \frac{\ln(1-x^2)}{x} \,\mathrm{d} x $$

\emph{[Changing the variable of integration does not alter the value of a definite integral]}

$$J=2\int_0^1 \frac{\ln(1-x)}{x}  \,\mathrm{d} x + 2\int_0^1 \frac{\ln(1+x)}{x}  \,\mathrm{d} x $$
$$J=2J+2\int_0^1 \frac{\ln(1+x)}{x}  \,\mathrm{d} x $$
$$2\int_0^1 \frac{\ln(1+x)}{x}  \,\mathrm{d} x =-J $$

From $(3)$,
\begin{equation}
I(-2)=-J
\end{equation}

Now,
$$I(\alpha)=\int_0^1 \frac{ln{(1- \alpha x+x^2)}}{x}  \,\mathrm{d} x$$

Since $I(x,\alpha) and \frac{\mathrm{d}I}{ \,\mathrm{d} \alpha}$ are continous in the interval $x \in (0,1)$ and are bounded, we can differentiate under the integral sign. Differentiating under the integral sign with respect to $\alpha$, we have:
$$I'(\alpha)=\int_0^1 \frac{1}{x} . \frac{-x}{(1- \alpha x+x^2)} \,\mathrm{d} x$$
$$I'(\alpha)=-\int_0^1 \frac{1}{(1- \alpha x+x^2)}   \,\mathrm{d} x$$
$$I'(\alpha)=-\int_0^1 \frac{1}{{(x-\frac{\alpha}{2})}^2 + \sqrt{(1-\frac{\alpha ^2}{4})}}   \,\mathrm{d} x$$

On evaluating the integral, we have:
$$I'(\alpha)=\frac{2}{\sqrt{4- \alpha ^2}} \left( \arctan(- \alpha)-\arctan\sqrt\frac{2- \alpha}{2+ \alpha}  \right) $$
$$ \,\mathrm{d} I(\alpha)=\left(\frac{2}{\sqrt{4- \alpha ^2}} \left( \arctan(- \alpha)-\arctan\sqrt\frac{2-\alpha}{2+ \alpha }  \right)\right)  \,\mathrm{d}  \alpha $$

Integrate both sides of the equation with respect to $\alpha$, from $\alpha=-2$ to $\alpha=2.$

Since $\frac{2}{\sqrt{4-\alpha^2}} \arctan(-\alpha)$ is an odd function, its integral over this interval will

be $0$.

$$I(2)-I(-2)=-2\int_{-2}^{2} \frac{1}{\sqrt{4- \alpha ^2}}\arctan\sqrt\frac{2-\alpha}{2+\alpha} \, \mathrm{d} \alpha$$

From (4) and (5),
$$I(2)-I(-2)=3J$$

Thus,
$$3J=-2\int_{-2}^{2} \frac{1}{\sqrt{4-\alpha ^2}}\arctan\sqrt\frac{2-\alpha}{2+\alpha} \,\mathrm{d} \alpha$$

Let $$u=\arctan\sqrt\frac{2- \alpha }{2+ \alpha }$$

Thus,
$$\, \mathrm{d} u=\frac{-1}{2\sqrt{4- \alpha ^2}} \, \mathrm{d} \alpha $$

When $\alpha \rightarrow -2$, $u \rightarrow \frac{\pi}{2}$ and when $\alpha \rightarrow 2$, $u\rightarrow 0$.
$$3J=4\int_{\frac{\pi}{2}}^0 u  \,\mathrm{d} u\ $$
$$3J=-\frac{\pi^2}{2}$$
$$J=-\frac{\pi^2}{6}$$

But, $J=-f(1).$
Thus,
$$f(1)=\frac{\pi^2}{6}$$

Since, $$f(x)=\sum_{n=1}^{\infty} {\frac{x^n}{n^2}} $$

it follows that

$$ f(1)= \sum_{n=1}^{\infty} {\frac{1}{n^2}} =\zeta(2) $$

Thus,
$$ \zeta(2)=\frac{\pi^2}{6}$$

\section{An Elementary Method To Calculate Values Of $\zeta(s)$ And Other Related Functions}

In this section, an elementary method is used to calculate the aforementioned values. An appropriate function is chosen and the following values are obtained for $s\in N$:

When $s$ is an even natural number, $$\zeta(s)= \sum_{n=1}^{\infty} {\frac{1}{n^s}}\ $$

and,

When $s$ is an odd natural number, let $$\Phi(s)= \sum_{n=1}^{\infty} {\frac{{(-1)}^n}{{(2n+1)}^s}}\ $$
Now, consider $\zeta(s)$ for $s \in N$, $s \geq 2$,
$$\zeta(s)=\sum_{n=1}^{\infty} \frac{1}{n^s} $$

Since, $n$ can be either odd or even,
$$\zeta(s)=\sum_{n=1}^{\infty} \frac{1}{{(2n)}^s} + \sum_{n=0}^{\infty} \frac{1}{{(2n+1)}^s} $$
$$\zeta(s)=\frac{1}{2^s}\sum_{n=1}^{\infty} \frac{1}{n^s} +  \sum_{n=0}^{\infty} \frac{1}{{(2n+1)}^s} $$
$$\zeta(s)=2^{-s} \zeta(s) + \sum_{n=0}^{\infty} \frac{1}{{(2n+1)}^s} $$

Thus,
\begin{equation}
\zeta(s)=\frac{2^s}{2^s-1}  \sum_{n=0}^{\infty} \frac{1}{{(2n+1)}^s}
\end{equation}

The Fourier Series of a periodic function $f(x)$, integrable on $[-\pi,\pi]$ is given:
$$f(x)=\frac{a_0}{2} + \sum_{n=1}^{\infty} {a_n \cos(nx)} + \sum_{n=1}^{\infty} {b_n \sin(nx)}\ $$

where,

for $n \geq 0$ $$a_n=\frac{1}{\pi} \int_{-\pi}^{\pi} f(x) \cos(nx)  \,\mathrm{d} x $$

and,
for $n \geq 1$
$$b_n=\frac{1}{\pi} \int_{-\pi}^{\pi} f(x) \sin(nx)  \,\mathrm{d} x$$

Consider the following definition of $f(x)$ which is periodic in $2\pi$:

$$
f(x) = \left\{
        \begin{array}{ll}
            |x| & \quad {-\pi} \leq x \leq {\pi} \\
            f(x+2k\pi) & \quad x > \pi, x<{-\pi}
        \end{array}
    \right.
$$

where $k \in N.$

Thus,
$$a_0=\frac{1}{\pi} \int_{-\pi}^{\pi} |x|  \,\mathrm{d} x $$
$$a_0=\frac{1}{\pi} \left (\int_{\pi}^0 x  \,\mathrm{d} x + \int_{0}^{-\pi} -x  \,\mathrm{d} x \right) $$
$$a_0=\pi$$

And,

$$a_n=\frac{1}{\pi} \int_{-\pi}^{\pi} |x| \cos(nx)  \,\mathrm{d} x $$
$$a_n=\frac{1}{\pi} \left (\int_{\pi}^0 x \cos(nx)  \,\mathrm{d} x + \int_{0}^{-\pi} -x \cos(nx)  \,\mathrm{d} x \right) $$
$$a_n=\frac{2}{n^2\pi} \left[ n \pi \sin(n \pi) +\cos(n \pi) -1 \right]$$

Since $n \in N,$ $$\sin(n \pi)=0$$ 
$$\cos(n\pi)={(-1)}^n.$$

Thus,
$$a_n=\frac{2}{n^2 \pi} \left[ {(-1)}^n - 1 \right] $$

Also,
$$b_n=\frac{1}{\pi} \int_{-\pi}^{\pi} |x| \sin(nx)  \,\mathrm{d} x $$

Since $|x|\sin(nx)$ is an odd function, its integral over $[-\pi,\pi]=0.$
Thus,
$$b_n=0$$

Therefore, the Fourier Series expansion for $f(x)$ defined above is:
$$|x|=\frac{a_0}{2} + \sum_{n=1}^{\infty} {a_n \cos(nx)} + \sum_{n=1}^{\infty} {b_n \sin(nx)}\ $$
$$|x|=\frac{\pi}{2} + \sum_{n=1}^{\infty} \frac{2}{n^2 \pi} \left[ {(-1)}^n - 1\right]\cos(nx)$$

When $n$ is even, $$\frac{2}{n^2 \pi} \left[ {(-1)}^n - 1\right]=0$$

When $n$ is odd, $$\frac{2}{n^2 \pi} \left[ {(-1)}^n - 1\right]=\frac{-4}{n^2}$$

Thus, $\forall x \in [-\pi,\pi]$,

$$ |x|=\frac{\pi}{2} - \frac{4}{\pi}\sum_{n=0}^{\infty} \frac{\cos{[(2n+1)x]}}{{(2n+1)}^2} $$

In what follows, we shall consider $x \geq 0$ only.

$\forall x \geq 0, |x|=x. $

Thus, $\forall x \in [0,\pi]$,

$$ x=\frac{\pi}{2} - \frac{4}{\pi}\sum_{n=0}^{\infty} \frac{\cos{[(2n+1)x]}}{{(2n+1)}^2} $$

Rearranging the terms of the equation, we have,
\begin{equation}
\sum_{n=0}^{\infty} \frac{\cos{[(2n+1)x]}}{{(2n+1)}^2} = \frac{\pi^2}{8}-\frac{\pi x}{4}
\end{equation}

Integrate both sides of this equation, with respect to $x$, from $x=0$ to $x=a$.
$$ \int_{0}^a \sum_{n=0}^{\infty} \frac{\cos{[(2n+1)x]}}{{(2n+1)}^2}  \,\mathrm{d} x = \int_{0}^a \left( \frac{\pi^2}{8}-\frac{\pi x}{4} \right)  \,\mathrm{d} x $$
$$ \sum_{n=0}^{\infty} \frac{\int_{0}^a \cos{[(2n+1)x]}  \,\mathrm{d} x} {{(2n+1)}^2}= \frac{\pi^2 a}{8} - \frac{\pi a^2}{8} $$
\begin{equation}
\sum_{n=0}^{\infty} \frac{\sin{[(2n+1)a]}}{{(2n+1)}^3} = \frac{\pi^2 a}{8} - \frac{\pi a^2}{8}
\end{equation}

Integrate both sides of this equation, with respect to $a$, from $a=x$ to $a=\frac{\pi}{2}$.
$$ \int_{x}^{\frac{\pi}{2}} \sum_{n=0}^{\infty} \frac{\sin{[(2n+1)a]}}{{(2n+1)}^3}  \,\mathrm{d} a = \int_{x}^{\frac{\pi}{2}} \left(\frac{\pi^2 a}{8} - \frac{\pi a^2}{8}\right)  \,\mathrm{d} a$$
$$ \sum_{n=0}^{\infty} \frac{\int_{x}^{\frac{\pi}{2}} \sin{[(2n+1)a]}  \,\mathrm{d} a}{{(2n+1)}^3}= \frac{\pi^4}{96} - \frac{\pi^2 x^2}{16} + \frac{\pi x^3}{24} $$
\begin{equation}
\sum_{n=0}^{\infty} \frac{\cos{[(2n+1)x]}}{{(2n+1)}^4} = \frac{\pi^4}{96} - \frac{\pi^2 x^2}{16} + \frac{\pi x^3}{24}
\end{equation}

Integrate both sides of this equation, with respect to $x$, from $x=0$ to $x=a$.
$$ \int_{0}^a \sum_{n=0}^{\infty} \frac{\cos{[(2n+1)x]}}{{(2n+1)}^4}  \,\mathrm{d} x = \int_{0}^a \left(\frac{\pi^4}{96} - \frac{\pi^2 x^2}{16} + \frac{\pi x^3}{24}\right)  \,\mathrm{d} x $$
$$ \sum_{n=0}^{\infty} \frac{\int_{0}^a \cos{[(2n+1)x]}  \,\mathrm{d} x} {{(2n+1)}^4}= \frac{\pi^4 a}{96} - \frac{\pi^2 a^3}{48} + \frac{\pi a^4}{96} $$
\begin{equation}
\sum_{n=0}^{\infty} \frac{\sin{[(2n+1)a]}}{{(2n+1)}^5} = \frac{\pi^4 a}{96} - \frac{\pi^2 a^3}{48} + \frac{\pi a^4}{96} 
\end{equation}

Integrate both sides of this equation, with respect to $a$, from $a=x$ to $a=\frac{\pi}{2}$.
$$ \int_{x}^{\frac{\pi}{2}} \sum_{n=0}^{\infty} \frac{\sin{[(2n+1)a]}}{{(2n+1)}^5}  \,\mathrm{d} a = \int_{x}^{\frac{\pi}{2}} \left( \frac{\pi^4 a}{96} - \frac{\pi^2 a^3}{48} + \frac{\pi a^4}{96} \right)  \,\mathrm{d} a $$
$$ \sum_{n=0}^{\infty} \frac{\int_{x}^{\frac{\pi}{2}} \sin{[(2n+1)a]}  \,\mathrm{d} a}{{(2n+1)}^5}= \frac{\pi^6}{960}-\frac{\pi x^5}{480}+\frac{\pi^2 x^4}{192} - \frac{\pi^4 x^2}{192} $$
\begin{equation}
\sum_{n=0}^{\infty} \frac{\cos{[(2n+1)x]}}{{(2n+1)}^6} = \frac{\pi^6}{960}-\frac{\pi x^5}{480}+\frac{\pi^2 x^4}{192} - \frac{\pi^4 x^2}{192}
\end{equation}

Integrate both sides of this equation, with respect to $x$, from $x=0$ to $x=a$.
$$ \int_{0}^a \sum_{n=0}^{\infty} \frac{\cos{[(2n+1)x]}}{{(2n+1)}^6}  \,\mathrm{d} x = \int_{0}^a \left(\frac{\pi^6}{960}-\frac{\pi x^5}{480}+\frac{\pi^2 x^4}{192} - \frac{\pi^4 x^2}{192}\right)  \,\mathrm{d} x $$
$$ \sum_{n=0}^{\infty} \frac{\int_{0}^a \cos{[(2n+1)x]}  \,\mathrm{d} x} {{(2n+1)}^6}= \frac{\pi^6 a}{960} - \frac{\pi x^6}{2880} + \frac{\pi^2 x^5}{960} - \frac{\pi^4  x^3}{576} $$
\begin{equation}
\sum_{n=0}^{\infty} \frac{\sin{[(2n+1)a]}}{{(2n+1)}^7} = \frac{\pi^6 a}{960} - \frac{\pi a^6}{2880} + \frac{\pi^2 a^5}{960} - \frac{\pi^4  a^3}{576}
\end{equation}

In the calculations performed above, the integral and summation sign have been interchanged which can be justfied by using a coroallry of the monotone and dominated convergence theorem. The dominated convergence theorem implies that we may integrate the sequence of partial sums term-by-term, which is tantamount to saying that we may switch integration and summation.

In equation (7), put $x=0.$
$$ \sum_{n=0}^{\infty} \frac{1}{{(2n+1)}^2} = \frac{\pi^2}{8} $$

From (6),
$$\zeta(2)=\frac{4}{3}.\frac{ {\pi}^2 }{8}= \frac{\pi^2}{6}$$

In equation (9), put $x=0.$
$$ \sum_{n=0}^{\infty} \frac{1}{{(2n+1)}^4} = \frac{\pi^4}{96} $$

From (6),
$$\zeta(4)=\frac{16}{15}.\frac{\pi^4}{96}=\frac{\pi^4}{90} $$

In equation (11), put $x=0.$
$$ \sum_{n=0}^{\infty} \frac{1}{{(2n+1)}^6} = \frac{\pi^6}{960} $$

From (6),
$$\zeta(6)=\frac{64}{63}.\frac{\pi^6}{960}=\frac{\pi^6}{945} $$

In equation (8),(10) and (12), put $x=\frac{\pi}{2}.$
Since,
$$ \sin\left[{(2n+1)}\frac{\pi}{2}\right]={(-1)}^n, $$

We have,
$$\Phi(3)=\sum_{n=0}^{\infty} \frac {{(-1)}^n} {{(2n+1)}^3} = \frac{\pi^3}{32} $$

And,
$$\Phi(5)=\sum_{n=0}^{\infty} \frac {{(-1)}^n} {{(2n+1)}^5} = \frac{5 \pi^5}{1536} $$

And,
$$\Phi(7)=\sum_{n=0}^{\infty} \frac {{(-1)}^n} {{(2n+1)}^7} = \frac{61 \pi^7}{184320} $$

Thus, the values of $\zeta(2)$ , $\zeta(4)$ and $\zeta(6)$ along with
$\Phi(s)$ for $s=3,5,7$ have been derived.

\section{Conclusion}
Thus, elementary methods have been used to arrive at values of the Riemann Zeta function at even natural numbers and a closely related infinite sum for odd natural numbers. The first method involved the technique of ''Differentiation Under The Integral Sign'' while the second used a Fourier Series expansion to arrive at the needed values. The second method can be used to attain the needed values for other integers as well by continuing the same process, in the same pattern i.e. by integrating over an interval and substituting the required values. The second method is efficient in a way that you get the sum of two series using one procedure, and can be continued endlessly to attain values that one needs. It may even be programmed to obtain results for large integers.

\end{document}